\documentclass[12pt,reqno]{amsart}
\usepackage{amssymb}
\usepackage{amsxtra}
\usepackage[mathscr]{eucal}

\theoremstyle{plain}

\theoremstyle{definition}

\theoremstyle{remark}

\begin{document}

\title[David M.~Bradley]{An infinite series for the natural logarithm that converges
throughout its domain and 
makes concavity transparent}

\date{\today}

\author{David~M. Bradley}
\address{Department of Mathematics and Statistics\\
         University of Maine\\
         5752 Neville Hall
         Orono, Maine 04469-5752\\
         U.S.A.}
\email[]{bradley@math.umaine.edu, dbradley@member.ams.org}



\begin{abstract}
The natural logarithm can be represented by an infinite series that converges for all positive real
values of the variable, and which makes concavity patently obvious.  Concavity of the natural 
logarithm is known to imply, among other things, the fundamental inequality between the arithmetic and geometric mean.  
 \end{abstract}

\maketitle

\interdisplaylinepenalty=500

The inequality
   $\log x < x-1$, 
which holds for all positive real $x\ne 1$~\cite[Theorem 150, p.\ 106]{HLP} is a fundamental but often overlooked property of the natural logarithm.  In fact, it is equivalent to concavity.     
To see this, fix $a>0$, replace $x$ by $x/a$, and note that if $0<x\ne a$, then 
\[
   \log x - \log a = \log(x/a) < (x/a)-1 = (x-a)\log' a,
\]
which says that the graph of the natural logarithm lies below its tangent line at every point
$(a,\log a)$.  On the other hand, concavity of the natural logarithm is known to imply the
ubiquitous inequality between the arithmetic and geometric mean of a finite set of positive real numbers~\cite[p.\ 119]{HLP}.

The inequality 
$\log x<x-1$ for $0<x\ne 1$ is an immediate consequence of 
the integral representation for the natural logarithm, via the identity
\[
   x-1-\log x = \int_1^x\int_1^t u^{-2}\,du\,dt,
\]
since the double integral is obviously positive for all positive real $x\ne 1$.  Nevertheless, it may
be of interest to derive an alternative representation which, rather than employ integrals, expresses
the difference $x-1-\log x$ as an infinite series each of whose terms is obviously positive for 
$0<x\ne 1$.  

To this end, we begin for $x>0$ with the definition of the derivative
\[
    \log x = \frac{d}{dh} x^h \bigg|_{h=0} = \lim_{h\to0}\frac{x^h-1}{h}
    = \lim_{h\to0} h^{-1}\big(x^h-1\big) = \lim_{n\to\infty}2^n\big(x^{2^{-n}}-1\big), 
\]
as a limit of a difference quotient~\cite[Ch.\ III, \S 6, pp.\ 175--176]{Courant}.  Using this, we can express $x-1-\log x$ for $x>0$ as a telescoping
series:
\begin{align*}
   x-1-\log x 
   & =\lim_{n\to\infty}\sum_{k=1}^n
   \big[2^{k-1}\big(x^{2^{1-k}}-1\big)-2^{k}\big(x^{2^{-k}}-1\big)\big]\\
   &=\sum_{k=1}^\infty\big[2^{k-1}\big(x^{2^{-k}}-1\big)\big(x^{2^{-k}}+1\big)-2^k
   \big(x^{2^{-k}}-1\big)\big]\\
   &=\sum_{k=1}^\infty 2^{k-1}\big(x^{2^{-k}}-1\big)\big[\big(x^{2^{-k}}+1\big)-2\big]\\
   &=\sum_{k=1}^\infty 2^{k-1}\big(x^{2^{-k}}-1\big)^2.
\end{align*}

The series converges by limit comparison with the convergent geometric series 
$\displaystyle\sum_{k=1}^\infty 2^{-k}$, for if $x>0$, then
\[
   \lim_{k\to\infty} 2^{k-1}\big(x^{2^{-k}}-1\big)^2 \bigg/ 2^{-k}
   = \lim_{k\to\infty} 2^{-k}2^{k-1}\bigg(\frac{x^{2^{-k}}-1}{2^{-k}}\bigg)^2
   = \frac12\big(\log x\big)^2.
\]
Since each term of the series is clearly positive for $0<x\ne 1$,  the 
inequality $\log x<x-1$
follows immediately.


\end{document}